\documentclass{amsart}
\usepackage[super]{cite}
\usepackage{amssymb,latexsym, dsfont, amsfonts, amsbsy, amsmath, mathrsfs, dsfont}    
\usepackage{tikz}

\usepackage{amsaddr}

\usepackage{thmtools}
\declaretheoremstyle[notefont=\bfseries,notebraces={}{},%
headpunct={},postheadspace=1em]{mystyle}
\declaretheorem[style=mystyle,numbered=no,name=Lemma]{lemma-hand}

\usepackage{thmtools}
\declaretheoremstyle[notefont=\bfseries,notebraces={}{},%
headpunct={},postheadspace=1em]{mystyle}
\declaretheorem[style=mystyle,numbered=no,name=Corollary]{cor-hand}

\theoremstyle{plain}
\newtheorem{theorem}{Theorem}[page]
\newtheorem{lemma}[theorem]{Lemma}
\newtheorem*{lemma*}{Lemma}
\newtheorem{corollary}[theorem]{Corollary}
\newtheorem{proposition}[theorem]{Proposition}

\theoremstyle{definition}

\newtheorem{definition}[theorem]{Definition}
\theoremstyle{remark}

\let\alg=\mathbf                        
\let\mod=\mathfrak							
	
\let\lang=\mathit

\begin{document}
\markboth{Dellunde and Vidal}{Truth-preservation under fuzzy pp-formulas}

\title{Truth-preservation under fuzzy pp-formulas}
%

\title{Truth-preservation under fuzzy pp-formulas}
\author{Pilar Dellunde}

\address{
	Philosophy Dpt. of the Autonomous University of Barcelona, \\
	Graduated School of Mathematics \& IIIA - CSIC\\
Bellaterra, Spain\\
	\footnote{pilar.dellunde@uab.cat}}

\author{Amanda Vidal}

\address{Institute of Computer Science, Czech Academy of Sciences\\Praha, Czech Republic\\
	\footnote{amanda@cs.cas.cz}}

%

\maketitle             

\begin{abstract}
How can non-classical logic contribute to the analysis of complexity in computer science? In this paper, we give a step towards this question, taking a logical model-theoretic approach to the analysis of complexity in fuzzy constraint satisfaction. We study fuzzy positive-primitive sentences, and we present an algebraic characterization of classes axiomatized by this kind of sentences in terms of homomorphisms and direct products. The ultimate goal is to study the expressiveness and reasoning mechanisms of non-classical languages, with respect to constraint satisfaction problems and, in general, in modelling decision scenarios.
\end{abstract}

\section{Introduction}
It is an natural question to observe the way in which mathematical non-classical logics can contribute to the study of certain phenomena, including complexity questions, in computer science. In particular, since many of the questions in CS are formulated over relational structures, the understanding of many-valued model theory, still in an early research stage, is very relevant for that purpose.  

If we focus in computational complexity questions, a first approach to the topic can be found in \cite{HorMoVi17} where some open problems are proposed by the authors about the relationship between fuzzy logic and valued constraint satisfaction. 
In our opinion, a research oriented to find a non-classical logical approach to complexity, should address, at least, the following three issues:

\begin{enumerate}
\item Show that there is a good trade-off between algebra and logic in the relevant fragments.
\item Identify which problems in complexity theory are naturally expressed as questions about the expressive power of the non-classical logic.
\item Prove that these complexity problems are not better addressed in other known logical formalisms. 
\end{enumerate}

These issues are naturally interrelated. To evaluate the trade-off between algebra and logic, it is important to identify which are the relevant fragments of the non-classical logic where the complexity problems have to be expressed; and to prove the relevancy of the fragments, a comparative study of different logical formalisms with respect to their expressive power has to be performed. 

Constraint-based modelling has become a central research area in computational social choice, and in particular in preference modelling, where preferences can be seen as soft constraints \cite{MeRoSchiex06} . Different soft constraint formalisms can be found in the literature, among which fuzzy constraint satisfaction (\cite{DuFarPra93} \cite{Rutt94},), possibilistic \cite{MouPra98} and probabilistic \cite{FarLang93} as well as a very general formulation over semirings \cite{BisMonRo97} that has not been, however, thoughtfully studied after. The most prominent case is that of Valued CSP (\cite{SchiexFarVer} \cite{KroZiv17},\cite{KolKroRo},), intensively studied with algebraic techniques and over which a plethora of complexity results have been proven.

The classical constraint satisfaction problem (CSP) is known to have strong connections with various problems in database theory and classical finite-model theory \cite{KoVar08} where CSP can be rephrased as the homomorphism problem, the conjunctive-query evaluation problem, or the join-evaluation problem (among others). Some problems in complexity theory are naturally expressed as questions about the expressive power of certain classical logics. For instance, Fagin's Theorem in Descriptive Complexity Theory (see \cite{Fa74}), shows that existential second-order logic captures NP on the class of all finite structures, that is, given an isomorphism-closed class K of finite structures of some fixed non-empty finite vocabulary, K is in NP if and only if K is definable by an existential second-order sentence.

When relaxation of constraints is considered, it is natural to ask for the relationship between the existing weighted CSP formalisms and the non-classical logical ones. In particular, as pointed out in \cite{HorMoVi17} the most usual weighted CSP formalisms (valued CSP and fuzzy CSP) can be equivalently formulated over predicate structures respectively of the product standard algebra and of G\"odel logic. This leads naturally to wonder the interests of model theory of predicate fuzzy logics in the understanding of relaxed constraint satisfaction questions. 
However, only in recent times, model theory of predicate fuzzy logics has been developed as a subarea of MFL (see for instance \cite{CiHa10} or \cite{DeGaNo16}), leaving the important area of fuzzy finite-model theory -the one nearest to the above questions- yet unexplored. Considering a general semantics for MFL, a plethora of left continuous t-norms can be defined, we can go far beyond of the minimum t-norm in the interval $[0,1]$ of the reals (most commonly used in FCSP) or the bounded addition in the negative natural numbers with bottom (arising in VCSP). We would like rather to explore the logical properties of fuzzy languages in general, their expressiveness, and reasoning mechanisms with respect to constraint satisfaction problems. \bigskip

Positive-primitive formulas are one of the key elements in the logical study of classical CSP (see for instance \cite{KoVar08}). In model theory applied to algebra, they have been also extensively studied. Algebraically, a pp-formula expresses the solvability of a system of linear equations, by asking for the existence of an assignment satisfying all atomic formulas conforming the pp-sentence. Particularly relevant is the use of  pp-formulas in model theory of modules, where every definable subset of a module is a boolean combination of pp-definable cosets. This fact is used to prove that the theory of modules has pp-elimination of quantifiers, introducing a useful logical tool to study these mathematical structures (for a reference see \cite{Prest88}).

In CSP, it is very important the relation of constraints (expressed with pp-formulas) and the clone of polymorphisms of the structure, which has led, among many other results, to conclude a remarkable complexity dichotomy of the CSP question for a given structure (\cite{Bul17} \cite{Zh17},). In the VCSP framework, the usual formulation of a constraint instance is that related to existential formulas of (strong conjunction) of atomic formulas \cite{KolKroRo} \cite{KroZiv17}, and a parallel relation with the so-called fractional polymorphisms is developed towards the complexity study of said problems. On the other hand, in \cite{HorMoVi17} a larger class of problems is treated (over certain algebraic classes), which are those related to existential formulas with both strong and weak conjunctions. 
 
The present paper revises and extends the results presented at the MDAI'18 Conference (see \cite{De18}). The original contribution of the article is the mathematical proof of two axiomatization theorems (for primitive-positive theories, and for existential positive theories). 

\section{Preliminaries}

\subsection{Predicate Fuzzy Logics} 
Now we present the syntax and semantics of the predicate fuzzy logic MTL$\forall_=$, one of the predicate extensions of the left-continuous t-norm based logic MTL introduced in \cite{EsGo01} and we refer to [Section 5, Chapter 1]\cite{CiHaNo11} for a complete and extensive presentation of MTL$\forall_=$. 
\begin{definition} \textbf{\textit{Syntax of Predicate Languages}}
 A \emph{predicate language} $\lang{P}$ is a triple $\left\langle Pred_{\lang{P}},Func_{\lang{P}},Ar_{\lang{P}} \right\rangle$, where $Pred_{\lang{P}}$ is a nonempty set of \emph{predicate symbols}, $Func_{\lang{P}}$ is a set of \emph{function symbols} (disjoint from $Pred_{\lang{P}}$), and $Ar_{\lang{P}}$ represents the \emph{arity function}, which assigns a natural number to each predicate symbol or function symbol. We call this natural number the \emph{arity of the symbol}. The predicate symbols with arity zero are called \emph{truth constants}, while the function symbols whose arity is zero are named \emph{individual constants}. 
\end{definition}

The set of $\lang{P}$-terms, $\lang{P}$-formulas and the notions of free occurrence of a variable, open formula, substitutability and sentence are defined as in classical predicate logic. We asume that the equality symbol $\approx$ of the language is interpreted in every structure as the crisp identity. 
Notice that, in the language we have introduced there are also function symbols. The results we present in this paper hold also for arbritrary languages, and for this reason we have presented a general proof, that could be used in further applications of pp-definability in non-relational structures, not necessarily related to Constraint Satisfaction Problems.

\begin{definition} We introduce an axiomatic system for the predicate logic \emph{MTL}$\forall_=$ 
\begin{description}
 
 \item[($\boldsymbol{\mathrm{P}}$)]$\space\space$ $\space\space$ $\space\space$ Instances of the axioms of the propositional logic \emph{MTL}.
 
 \item[($\boldsymbol{\forall 1}$)]$\space\space$ $(\forall x)\varphi(x)\rightarrow\varphi(t)$, where the term $t$ is substitutable for $x$ in $\varphi$.
 
 \item[($\boldsymbol{\exists1}$)]$\space\space$ $\varphi(t)\rightarrow(\exists x)\varphi(x)$, where the term $t$ is substitutable for $x$ in $\varphi$.
 
 \item[($\boldsymbol{\forall 2}$)]$\space\space$ $(\forall x)(\xi\rightarrow\varphi)\rightarrow(\xi\rightarrow(\forall x)\varphi(x))$, where $x$ is not free in $\xi$.
 
 \item[($\boldsymbol{\exists2}$)]$\space\space$ $(\forall x)(\varphi\rightarrow\xi)\rightarrow((\exists x)\varphi\rightarrow\xi)$, where $x$ is not free in $\xi$.
 
 \item[($\boldsymbol{\forall3}$)]$\space\space$ $(\forall x)(\varphi\vee\xi)\rightarrow(\varphi \vee (\forall x)\xi)$, where $x$ is not free in $\varphi$.

 \item[(Ref$_=$)]$\space\space$ $(\forall x)x=x$
  \item[(LP)]$\space\space$$(\forall x)(\forall y )(x=y \to (\psi(x) \to \psi(y))$
  \item[(CRISP$_=$)]$\space\space$$(\forall x)(\forall y )(x=y \vee \neg x=y)$
\end{description}

\bigskip

The deduction rules of \emph{MTL}$\forall_=$ are those of \emph{MTL} and the rule of generalization: from $\varphi$ infer $(\forall x)\varphi$. The definitions of proof and provability are analogous to the classical ones. A set of formulas $\Phi$ is \emph{consistent}, if $\Phi\not\vdash\overline{0}$.
 \end{definition}

\begin{definition} \label{evaluation}\textbf{\textit{Semantics of Predicate Fuzzy Logics}} Consider a predicate language $\lang{P}=\langle Pred_{\lang{P}}, Func_{\lang{P}}, Ar_{\lang{P}} \rangle$ and $\textbf{A}$ an MTL-algebra. 
	We define an $\alg{A}$\emph{-structure} 
	$\mod{M}$ for $\lang{P}$ as a triple 
	$\langle M, (P_{\mod{M}})_{P\in Pred}, (F_{\mod{M}})_{F\in Func} \rangle$, where $M$ is a nonempty domain, $P_{\mod{M}}$ is an $n$-ary fuzzy relation for each $n$-ary predicate symbol, identified with an element of $\textbf{A}$, if $n=0$; and $F_{\mod{M}}$ is a function from $M^n$ to $M$, identified with an element of $M$, if $n=0$. 

As usual, if $\mod{M}$ is an $\alg{A}$-structure for $\lang{P}$, an $\mod{M}$-evaluation of the object variables is a mapping $v$ assigning to each object variable an element of $M$. The set of all object variables is denoted by $Var$. If $v$ is an $\mod{M}$-evaluation, $x\in Var$ and $a\in M$, we denote by $v[x\mapsto a]$ the $\mod{M}$-evaluation so that $v[x\mapsto a](x)=a$ and $v[x\mapsto a](y)=v(y)$ for $y$ an object variable such that $y\not=x$. If $\mod{M}$ is an  $\alg{A}$-structure and $v$ is an $\mod{M}$-evaluation, we define the \emph{values} of terms, and the \emph{truth values} of formulas in $\mod{M}$ for an evaluation $v$ recursively as follows:
 
\begin{itemize}

\item $||x||_{\mod{M},v}=v(x)$;

\item $||F(t_1,\ldots,t_n)||_{\mod{M},v}=F_{\mod{M}}(||t_1||_{\mod{M},v},\ldots,||t_n||_{\mod{M},v})$, for $F\in Func$;

\item $||P(t_1,\ldots,t_n)||_{\mod{M},v}=P_{\mod{M}}(||t_1||_{\mod{M},v},\ldots,||t_n||_{\mod{M},v})$, for $P\in Pred$;

\item $||\lambda(\varphi_1,\ldots,\varphi_n)||_{\mod{M},v}=\lambda_{\alg{A}}(||\varphi_1||_{\mod{M},v},\ldots,||\varphi_n||_{\mod{M},v})$, for every connective $\lambda$;

\item $||(\forall x)\varphi||_{\mod{M},v}=inf\{||\varphi||_{\mod{M},v[x\rightarrow a]}\mid a\in M\}$;

\item $||(\exists x)\varphi||_{\mod{M},v}=sup\{||\varphi||_{\mod{M},v[x\rightarrow a]}\mid a\in M\}$.

\end{itemize}
 \end{definition} 

If the infimum or supremum does not exist, we take the truth value of the quantified formula as undefined. We say that the $\alg{A}$-structure $\mod{M}$ is \emph{safe} if $||\varphi||_{\mod{M},v}$ is defined for every formula and every evaluation. We restrict only to safe models. 

\bigskip

We assume that the language has an equality symbol $\approx$, interpreted as a crisp identity. We denote by $||\varphi||_{\mod{M}}=1$ the fact that $||\varphi||_{\mod{M},v}=1$ for all $\mod{M}$-evaluation $v$; and given a set of sentences $\Phi$, we say that $\mod{M}$ is a \emph{model of $\Phi$}, if for every $\varphi\in\Phi$, 
$||\varphi||_{\mod{M}}=1$. A theory is a set of sentences, and we denote by $Mod _\alg{A}$($\Phi$) the set of $\alg{A}$-models of $\Phi$, and by Th($\mod{M}$), the theory of $\mod{M}$, that is, the set of sentences evaluated 1 in $\mod{M}$. We say that two models are \emph{elementary equivalent}, if they have the same theory.

\bigskip

From now on we fix a finite linearly ordered MTL-algebra \textbf{A} and consider only structures over this algebra. Examples of this type of finite algebras are all the \L ukasiewicz finite-valued algebras. Since we work with structures over a fixed finite linearly ordered MTL-algebra, the infimum and the supremum in Definition \ref{evaluation} always exist, and they coincide with the minimum and maximum. Not only this, but the linear ordering of the underlying algebra makes this values coincide with one of the elements in consideration. 
This is called \textit{witnessing condition}, and in particular for the development of this article, we are interested in the $\exists$-witnessing: for every formula of the form $(\exists \overline{x}) \psi (\overline{x})$, there are $\overline{d} \in M$ such that $||(\exists \overline{x}) \psi  (\overline{x})||_{\mod{M}}=||\psi  (\overline{d})||_{\mod{M}}$. The witnessing condition, not necessarily holding neither over non-linearly ordered algebras nor over infinite ones\footnote{Though there are some particular cases of logics whose models are not necessarily witnessed, but they enjoy completeness with respect to witnessed models of the class \cite{Ha07, Ha07a}.}, is an important requirement in the results presented below. It allows to preserve truth under transformations that are only determining the behaviour of atomic formulas over the top element of the algebra, while in the case of not witnessed structures, that is not specific enough to conclude analogous results.

\bigskip
Another important property when we work with structures over a finite linearly ordered MTL-algebras is compactness, both for satisfiability and consequence (the proof can be found in [Th. 4.4]\cite{De14}). Remark that, in fuzzy logic it is not always the case, for instance the product predicate logic is neither satisfiability nor consequence compact with respect to its standard algebra. Given a set of sentences $\Sigma$, and a sentence $\phi$, we denote by $\Sigma \models_{\mathbf{A}}\phi$ the fact that every \textbf{A}-model of $\Sigma$ is also an \textbf{A}-model of $\phi$.

\bigskip

\begin{theorem}\label{com}\textbf{\textit{\textbf{A}-compactness.}}
	 For every set of sentences $\Sigma$ and sentence $\phi$, the following holds:
\begin{enumerate} 
\item \emph{[Satisfiability]} If for every finite subset $\Sigma_0\subseteq\Sigma$, $\Sigma_0$ has an \emph{\textbf{A}}-model, then $\Sigma$ has also an \emph{\textbf{A}}-model.
\item  \emph{[Consequence]} If $\Sigma \models_{\mathbf{A}}\phi$, then there is a finite subset $\Sigma_0\subseteq\Sigma$ such that $\Sigma_0 \models_{\mathbf{A}}\phi$.
\end{enumerate}
\end{theorem}

\begin{definition} \label{st} Let $\lang{P}$ be a predicate language and let $\lang{P'}$ an expansion of $\lang{P}$ that contains constant symbols $d_a$ for some elements of the algebra $a \in A$. We say that a structure $\mod{M}$ for the expanded language $\lang{P'}$ is \emph{standard} if every algebra constant symbol is interpreted as the corresponding element of the algebra.
\end{definition}

In the present paper we will make an extensive use of modelo-theoretic techniques (such as the diagrams method, for instance) that work with expanded languages. The compactness result of [Th. 4.4]\cite{De14} hold also for languages expanded adding new algebra constant symbols for the elements of \textbf{A}. The important fact we want to point out here is that, if every finite subset $\Sigma_0\subseteq\Sigma$ has an standard model, we can guarantee that, using \textbf{A}-compactness, $\Sigma$ will also have an standard model. This is because the proof of [Th. 4.4]\cite{De14} uses an ultraproduct construction, and by the same definition of ultraproduct [Def. 3.2] \cite{De14} if we have a set of structures where every algebra constant symbol $\overline{a}$ is interpreted as the element of the algebra $a$, then the ultraproduct of this set of structures also interprets every $\overline{a}$ as the corresponding element of the algebra $a$. Otherwise stated, from now on we assume that all the structrures for the expanded languages with algebra constant symbols are standard (in the sense of Definition \ref{st}). Also, for the sake of simplicity, when it is clear from the context, we will refer to \textbf{A}-structures simply as \emph{structures}, because all the structures we consider will be over the same algebra \textbf{A}.

\section{Fuzzy Positive-Primitive Formulas} 

In (classical) CSP, pp-formulas understood as an existential positive (with only $\land$-operation, see \cite{h} for a general reference of the classical positive-primitive fragment) faithfully represent the CSP instances, since these are usually understood simply as a set of atomic formulas (constraints) and a solution for it is a tuple of elements that satisfy all the constraints at the same time. On the other hand, when relaxing this question to a valued setting, there is more than one natural definition of an instance of a valued constraint problem. While the most general definition of a VCSP instance (see eg. \cite{KroZiv17}) is analogous to that of CSP (a set of constraints, that now will be fuzzy and whose weights will be combined in some uniform way), we can also consider combinations of both conjunctions (as in \cite{HorMoVi17}), to combine the constraint weights in two different levels. Along this section we will introduce and study what we consider to be the three natural generalisations of pp-formulas to the valued context. In particular, we see that homomorphisms and direct products preserve all these positive-primitive formulas. 

\begin{definition}\label{d1:pp} \textbf{\textit{Fuzzy Positive-Primitive Formula}}
	Given a predicate language $\mathcal{P}$, and a $\mathcal{P}$-formula $\phi$, it is said that $\phi$ is 
	\begin{enumerate}
		\item \emph{$\wedge$-primitive formula} if $\phi$ is of the form $(\exists \overline{x})\psi$, where $\psi$ is a quantifier-free formula  built from atomic formulas by using only the connective $\wedge$.
		\item \emph{$\&$-primitive formula} if $\phi$ is of the form $(\exists \overline{x})\psi$, where $\psi$ is a quantifier-free formula  built from atomic formulas by using only the connective $\&$.
		\item \emph{$\wedge\&$-primitive formula} if $\phi$ is of the form $(\exists \overline{x})\psi$, where $\psi$ is a quantifier-free formula  built from atomic formulas by using only the connectives $\wedge$ and $\&$.
	\end{enumerate}
%
%
\end{definition} 
Observe that clearly $\wedge$-primitive formulas and $\&$-primitive formulas are both $\wedge\&$-primitive formulas too, so in order to simplify the reading we will sometimes refer to  $\wedge\&$-primitive formulas by pp-formulas.

%


 \bigskip
Using the following proposition of \cite{EsGo01} it can be proved that there is a canonical normal form for $\wedge\&$-primitive formulas. That is, every $\wedge\&$-primitive formula  is equivalent to one of the form $\exists x_1 \cdots x_n \bigwedge_{i \in I} \bigodot _{j \in J} \phi_{ij}(x_1 \ldots x_n)$, where for every $i \in I, j \in J$, $\phi_{ij}$ is an atomic formula.

\begin{proposition} \label{godo}\textbf{\textit{[Prop. 1.30 \cite{EsGo01}]}} Let $\mathcal{P}$ be a predicate language, for every $\mathcal{P}$-formulas $\phi, \psi$ and $\alpha$, the formula $\phi \& (\psi \wedge \alpha) \leftrightarrow ((\phi \& \psi ) \wedge (\phi \& \alpha))$ is an \emph{MTL}$\forall_=$ theorem.
\end{proposition}

 \bigskip

Let us recall now the definition of homomorphism introduced in \cite{DeGaNo16} as a generalization of the notion of classical homomorphism.
\begin{definition}\label{def:homo}\textbf{\textit{Homomorphism}} Let $\mathcal{P}$ be a predicate language, $\mod{M}$ and $\mod{N}$ be two structures for $\mathcal{P}$, and $g$ a mapping from $M$ to $N$. We say that $g$ is a \emph{homomorphism} from $\mod{M}$ to $\mod{M}$ if and only if
\begin{enumerate}
\item For every $n$-ary function symbol $F\in\mathcal{P}$, and $d_1,\ldots,d_n\in M$, $$g(F_{\mod{M}}(d_1,\ldots,d_n))=F_{\mod{N}}(g(d_1),\ldots,g(d_n)).$$
\item For every $n$-ary predicate symbol $P\in\mathcal{P}$, and $d_1,\ldots,d_n\in M$, 
\begin{center} if $||P(d_1,\ldots,d_n)||_{\mod{M}}=1$, then $||P(g(d_1),\ldots,g(d_n))||_{\mod{N}}=1$.
\end{center}
\end{enumerate}
Moreover, we say that $g$ is an \emph{embedding}, if $g$ is one-to-one, and that $g$ is an \emph{isomorphism}, if $g$ is a surjective embedding.
\end{definition}

In the following lemma we prove that primitive positive formulas are preserved under homomorphisms.

\begin{lemma} \label{lem:pp preserving} Let $\mathcal{P}$ be a predicate language, $\mod{M}$ and $\mod{N}$ be two structures for $\mathcal{P}$, $g$ a homomorphism from $\mod{M}$ to $\mod{M}$, and $\phi$ a pp-formula. Then, for every $d_1,\ldots,d_n\in M$, 
\begin{center} if $||\phi(d_1,\ldots,d_n)||_{\mod{M}}=1$, then $||\phi(g(d_1),\ldots,g(d_n))||_{\mod{N}}=1$.
\end{center}
\end{lemma}
 \begin{proof} By induction on the complexity of $\phi$. 
\smallskip

\noindent {\bf Atomic step.} Let $\phi$ be an atomic formula of the form $P(t_1 \ldots, t_k)$, where $P\in\mathcal{P}$ is a predicate symbol, and $t_1 \ldots, t_k$ are $\mathcal{P}$-terms. Since $g$ is a homomorphism, we have that, in general, for every $\mathcal{P}$-term $t$, and $d_1,\ldots,d_n\in M$, $g(t_{\textbf{M}}(d_1,\ldots,d_n))=t_{\textbf{N}}(g(d_1),\ldots, g(d_n))$ and thus

\medskip$\begin{array}{lr} ||P(t_1 \ldots, t_k)(d_1,\ldots,d_n)||_{\mod{M}}=1\Rightarrow
 
\\[1ex]  ||P(t_{1\textbf{M}}(d_1,\ldots,d_n), \ldots, t_{k\textbf{M}}(d_1,\ldots,d_n))||_{\mod{M}}=1\Rightarrow

\\[1ex]  ||P(g(t_{1\textbf{M}}(d_1,\ldots,d_n)), \ldots, g(t_{k\textbf{M}}(d_1,\ldots,d_n)))||_{\mod{N}}=1\Rightarrow

\\[1ex]||P(t_{1\textbf{N}}(g(d_1),\ldots,g(d_n)), \ldots, t_{k\textbf{N}}(g(d_1),\ldots,g(d_n))||_{\mod{N}}=1\Rightarrow

\\[1ex] ||P(t_1 \ldots, t_k)(g(d_1),\ldots,g(d_n))||_{\mod{N}}=1. & \end{array}$
  
\bigskip

\noindent {\bf Quantifier-free.} Assume inductively that the property holds for $\psi$ and for $\chi$, then we have: 

\medskip$\begin{array}{lr} 1=||\psi \& \chi (d_1,\ldots,d_n)||_{\mod{M}}=||\psi  (d_1,\ldots,d_n)||_{\mod{M}} * ||\chi (d_1,\ldots,d_n)||_{\mod{M}}\Rightarrow
 \\[1ex]  ||\psi  (d_1,\ldots,d_n)||_{\mod{M}} =1 \mbox{ and } ||\chi (d_1,\ldots,d_n)||_{\mod{M}}=1\Rightarrow
\\[1ex]  ||\psi (g(d_1),\ldots,g(d_n))||_{\mod{N}}  =1 \mbox{ and } ||\chi (g(d_1),\ldots,g(d_n))||_{\mod{N}}=1\Rightarrow
 \\[1ex]||\psi (g(d_1),\ldots,g(d_n))||_{\mod{N}} * ||\chi (g(d_1),\ldots,g(d_n))||_{\mod{N}}=1\Rightarrow
\\[1ex] ||\psi \& \chi (g(d_1),\ldots,g(d_n))||_{\mod{N}}=1. & \end{array}$

\bigskip

\noindent Observe that the same argument holds for the weak conjunction $\wedge$.

\smallskip

\noindent {\bf Existential step.} Assume inductively that the property holds for $\psi(x)$. Since $\mod{M}$ is an $\exists$-witnessed structure, we have that for some $e \in M$,
$$||(\exists x) \psi  (x,d_1,\ldots,d_n)||_{\mod{M}}=||\psi  (e, d_1,\ldots,d_n)||_{\mod{M}} $$
Thus, if $||(\exists x) \psi  (x,d_1,\ldots,d_n)||_{\mod{M}}=1$, then $||\psi  (e, d_1,\ldots,d_n)||_{\mod{M}}=1$ and, by inductive hypothesis, $$1=||\psi (g(e),g(d_1),\ldots,g(d_n))||_{\mod{N}}\leq ||(\exists x) \psi  (x,g(d_1),\ldots,g(d_n))||_{\mod{N}}.$$
\end{proof}  

Observe also formulas closed under $\vee$ are preserved under the previous family of homomorphisms. However, that is not the case with other operations (namely $\forall$ and $\rightarrow$). While taking subjective homomorphisms allows to preservve universal formulas too, the implication is hardly to have a regular behaviour.


Let us now introduce the notion of direct product. Unlike other definitions introduced in the literature, for instance in \cite{HorMoVi17} we are interested in products resulting in structures over the same algebra \textbf{A}. 
\begin{definition}\label{def:product}\textbf{\textit{\textbf{A}-direct product}} Let $\mathcal{P}$ be a predicate language, $I$ a nonempty set, and for every $i \in I$, $\mod{M}_i$  is a structure for $\mathcal{P}$. The \emph{direct product} of the family $\{\mod{M}_i:i \in I\}$,
denoted by $\prod_{i\in I} \mod{M}_i$, is the structure that has as domain the usual classical direct product, the usual classical interpretation for constants and function symbols, and for every n-adic predicate symbol $ P \in\mathcal{P}$ and tuple of elements $\overline{d_1},\ldots, \overline{d_n}$ of $\prod_{i\in I}M_i$,
$$P_{\prod_{i\in I} \mod{M}_i}(\overline{d_1},\ldots, \overline{d_n})=\min \{ P_{\mathbf{M}_i}(\overline{d_1}(i),\ldots, \overline{d_n}(i)):i \in I\}$$
\end{definition}

Notice that the product is well-defined because the algebra is finite. Moreover, observe that, so defined, the $i$-projection of the direct product onto $\mod{M}_i$ is a homomorphism, and thus, by Lemma \ref{lem:pp preserving}, preserves pp-formulas. Focusing in truth-preservation only, a second kind of product-like construction can be done, obtaining a class of structures of which \textbf{A}-direct product are a particular case.

\begin{definition} \textit{\textbf{Weak \textbf{A}-direct product}} Let $\mathcal{P}$ be a predicate language, $I$ a nonempty set, and for every $i \in I$, $\mod{M}_i$  is a structure for $\mathcal{P}$. A \emph{weak \textbf{A}-direct product} of the family $\{\mod{M}_i:i \in I\}$
 is any structure $\mod{M}$ that has as domain the usual classical direct product, and the usual classical interpretation for constants and function symbols, and for every n-adic predicate symbol $ P \in\mathcal{P}$, and tuples of elements $\overline{d_1},\ldots, \overline{d_n}$ of $\prod_{i\in I}M_i$,
$$ P_{\mod{M}}(\overline{d_1},\ldots, \overline{d_n})=1 \text{ if and only if }P_{\mod{M}_i}(\overline{d_1}(i),\ldots, \overline{d_n}(i))=1, \text{ for every }i \in I. $$

We will denote by $w\Pi_{i \in I}\mod{M}_i$ the family of weak \textbf{A}-direct products of $\{\mod{M}_i\}_{i \in I}.$
	\end{definition}
	
Observe that given a family of models, its \textbf{A}-direct product (in the sense of Definition \ref{def:product}) is always a weak $\alg{A}$-directed product.

\begin{lemma} \label{lem:prod preserving} Let $\mathcal{P}$ be a predicate language, $I$ a nonempty set, and for every $i \in I$, $\mod{M}_i$  a $\mathcal{P}$-structure. Assume that $\phi$ is a positive-primitive $\mathcal{P}$-formula, and $\overline{d_1},\ldots, \overline{d_n}$ are tuples of elements of $\prod_{i\in I}M_i$. Then for every $\mod{M} \in w\Pi_{i \in I}\mod{M}_i$ the following holds: for every $i \in I$, 
$$||\phi(\overline{d_1}(i),\ldots, \overline{d_n}(i))||_{\mathbf{M}_i}=1 \text{ if and only if }||\phi(\overline{d_1},\ldots, \overline{d_n})||_{\mod(M)}=1. $$
\end{lemma}
\begin{proof}  
By induction on the complexity of $\phi$. The proof of the atomic and quantifier-free step is similar to the corresponding proof in Lemma \ref{lem:pp preserving}, by using the fact that for every $\mathcal{P}$-term $t$, $$t_{\mod(N)}(\overline{d_1},\ldots, \overline{d_n})=(t_{\mod{M_i}}(\overline{d_1}(i),\ldots, \overline{d_n}(i)):i \in I)$$ For the existential step, assume inductively that the property holds for $\psi(x)$. For every $i \in I$, $||(\exists x)\psi(x,\overline{d_1}(i),\ldots, \overline{d_n}(i))||_{\mod{M}_i}=1$ (because the structures are $\exists$-witnessed) if and only if for every $i\in I$, there is  $\overline{e}(i)\in M_i$ such that $$||\psi(\overline{e}(i),\overline{d_1}(i),\ldots, \overline{d_n}(i))||_{\mod{M}_i}=1$$ Then, by using the inductive hypothesis, this happens if and only if $$1=||\psi(\overline{e},\overline{d_1},\ldots, \overline{d_n})||_{\mod(N)}\leq ||(\exists x) \psi (x,\overline{d_1}(i),\ldots, \overline{d_n}(i))||_{\mod(N)}.$$
\end{proof}

As a particular case, the previous holds for the $\alg{A}$-directed products.

\section{Fuzzy Positive-Primitive Sets of Axioms}

Axiomatization theorems provide a correspondence between syntactic and semantic notions in logic. Diagrams are the building blocks that, glued with compactness, allow us to build extensions of structures, and prove these axiomatization theorems. Let us thus to introduce the method of diagrams in this fuzzy setting in order to characterize homomorphisms, and prove an equivalent condition to the preservation of pp-formulas between structures.
\begin{definition}\label{lfirst} Let $\mathcal{P}$ be a predicate language, and $\mod{M}$ a structure for $\mathcal{P}$. The expansion of the language $\mathcal{P}$ by adding an individual constant symbol $c_m$ for every $m\in M$, is denoted by $\mathcal{P}^{M}$; and the expansion of the structure $\mod{M}$ to $\mathcal{P}^{M}$ is denoted by $\mod{M}^\sharp$, where for every $m\in M$, $(c_m)_{\mod{M}^\sharp}=m$.
\end{definition} 
\begin{definition}\label{d:adiag} Let $\mathcal{P}$ be a predicate language. For every structure $\mod{M}$ for $\mathcal{P}$, we define \emph{Diag}$(\mod{M})$ as the set of atomic $\mathcal{P}^{M}$-sentences $\sigma$ such that $||\sigma||_{\mod{M}^\sharp}=1$.
\end{definition}
Following the same lines of the proof of [Prop. 32]\cite{De11} we can obtain this characterization of homomorphisms in terms of diagrams.
\begin{lemma}\label{l:ordiag}  Let $\mathcal{P}$ be a predicate language, $\mod{M}$ and $\mod{N}$ be two structures for $\mathcal{P}$. The following are equivalent: 
\begin{enumerate} 
\item There is an expansion of $\mod{N}$ that is a model of \emph{Diag}$(\mod{M})$.
\item There is a homomorphism $g:M \to N$ from $\mod{M}$ to $\mod{N}$.
\end{enumerate}
\end{lemma}

Notice that, since the Diag$(\mod{M})$ contains equalities but not inequalities, the obtained homomorphism does not need to be an embedding. Now we present a characterization in terms of extensions, of when two structures preserve pp-formulas. Given a structure $\mod{N}$ for a language $\mathcal{P}$, we denote by $\mod{N}_{\mathbf{A}}$, the expansion of $\mod{N}$ to a language obtained adding to $\mathcal{P}$ a new truth constant $d_a$ for every element $a$ of the algebra, and such that $d_a$ is interpreted as the corresponding element $a$.
 
\begin{proposition} \label{lem:atax} Let $\mathcal{P}$ be a predicate language, and $\mod{M}$ and $\mod{N}$ be two structures for $\mathcal{P}$. Then, every pp-sentence which is evaluated $1$ in $\mod{M}$, is also evaluated $1$ in $\mod{N}$ if and only if there is a structure $\mod{L}$ for $\mathcal{P}$, elementary equivalent to $\mod{N}$, and a homomorphism $g$ from $\mod{M}$ to $\mod{L}$.
\end{proposition}

\begin{proof} The direction from left to right is clear. Now assume that every pp-sentence which is evaluated $1$ in $\mod{M}$, is also evaluated $1$ in $\mod{N}$. First we show that Diag$(\mod{M})\cup$Th$(\mod{N_{\mathbf{A}}})$ has a model. We prove that for every finite subset $\{\sigma_1 \dots, \sigma_n\}$ of Diag$(\mod{M})$, $\{\sigma_1 \dots, \sigma_n\}\cup $Th$(\mod{N}_{\mathbf{A}})$ has a model. 
Let $c_{m_1}, \ldots, c_{m_k}$ be the object constants of the expanded language that occur in $\{\sigma_1 \dots, \sigma_n\}$. For every $1\leq i \leq n$, let $\sigma'_i$ be the formula obtained from $\sigma_i$ by substituting the constants
$c_{m_1}, \ldots, c_{m_k}$ by new variables $\overline{y}=y_{m_1}, \ldots, y_{m_k}$. Then we have that $||(\exists\overline{y})(\sigma'_1 \wedge \cdots \wedge \sigma_n(\overline{y}))||_{\mod{M}}=1$ and thus, since every pp-sentence which is evaluated $1$ in $\mod{M}$, is also evaluated $1$ in $\mod{N}$, and $(\exists\overline{y})(\sigma'_1 \wedge \cdots \wedge \sigma'_n(\overline{y}))$ is a pp-sentence, we obtain
$||(\exists\overline{y})(\sigma'_1 \wedge \cdots \wedge \sigma_n(\overline{y}))||_{\mod{N}}=1$.

Since $\mod{N}$ is a $\exists$-witnessed structure, we have a sequence of elements of N, $\overline{e}=e_{m_1}, \ldots, e_{m_k}$, such that $||\sigma'_1 \wedge \cdots \wedge \sigma_n(\overline{e})||_{\mod{N}}=1$. If we assign to every constant $c_{m_i}$ the corresponding element $e_{m_i}\in N$ we obtain an expansion of $\mod{N}_{\mathbf{A}}$ that satisfies $\{\sigma_1 \dots, \sigma_n\}\cup $Th$(\mod{N}_{\mathbf{A}})$.

By $\mathbf{A}$-compactness for satisfiability, there is a structure $\mod{L}$ for $\mathcal{P}$ that has an expansion which is a model of Diag$(\mod{M})\cup$Th$(\mod{N}_{\mathbf{A}})$. By Lemma \ref{l:ordiag}, there is a homomorphism $g$ from $\mod{M}$ to $\mod{L}$. Moreover, since $\mod{L}$ is a model of Th$(\mod{N}_{\mathbf{A}})$, $\mod{L}$ is elementary equivalent to $\mod{N}$. \end{proof}

Notice that the previous proposition is true also when we subtitute in its statement \emph{pp-sentences} by \emph{$\wedge$-primitive sentences} or by \emph{$\&$-primitive sentences}.

\bigskip
Now we prove an axiomatization theorem for theories closed under homomorphisms and weak direct products. Recall that a theory $T$ is \emph{closed under} an operator $O$, if the class of its $\alg{A}$-models, $Mod _\alg{A}$($T$), is closed under $O$, that is, if $O$ is an n-adic operator, and $\mod{M}_1, \ldots, \mod{M}_n$ are $\alg{A}$-models of $T$, then $O(\mod{M}_1, \ldots, \mod{M}_n)$ is also an $\alg{A}$-model of $T$. And it is said that a theory $T$ is \emph{axiomatized} by a set of sentences $\Sigma$, if $Mod _\alg{A}$($T$)$=Mod _\alg{A}$($\Sigma$). 
\begin{theorem}\label{t:axiom} Let $\mathcal{P}$ be a predicate language and $T$ be a consistent theory. Then, $T$ is closed under homomorphisms and direct products if and only if $T$ is axiomatized by a set of positive primitive sentences.
\end{theorem}
 
\begin{proof} $\Leftarrow$ is immediate. For the other direction, consider the sets of formulas $ppT = \{\varphi: \varphi$ is a pp-sentence and $T\models_{\mathbf{A}}\varphi \}$, and $\overline{ppT} = \{\varphi: \varphi $ is a pp-sentence and $T\not \models_{\mathbf{A}}\varphi\}$.  We aim to prove that $ppT\models_{\mathbf{A}} T$. Clearly $ppT \neq \emptyset$. Now, consider a model $\mod{N}$ of $ppT$. Assume first that $\overline{ppT}  \neq \emptyset$. Then for each $\chi \in \overline{ppT}$, let $\mod{M}_\chi$ be a model of $T$ but not a model of $\chi$, and let $\mod{M} = \Pi_{\chi \in \overline{ppT}} M_\chi$. 
	 
Since $T$ is closed under direct products, $\mod{M}$ is a model of $T$ too. From Proposition \ref{lem:atax} we know that for any pp-sentence $\psi$, $\psi$ is valid in $\mod{M}$ if and only if $\psi$ is valid in $\mod{M}_\chi$ for all $\chi \in \overline{ppT}$. Allow us to call $pp\mod{M}_\chi$ the set of pp-sentences valid in $\mod{M}_{\chi}$. It is easy to see that $\bigcap_{\chi \in \overline{ppT}} pp\mod{M}_\chi = ppT$. Indeed, $\supseteq$ is trivial, since all $\mod{M}_\chi$ are models of $T$. On the other hand, for any pp-sentence $\chi$ not in $ppT$, $\chi \in \overline{ppT}$, and thus, by definition, $\chi \not \in pp\mod{M}_\chi$.

Now, we know that the pp-sentences valid in $\mod{M}$ are exactly $ppT$. Thus, $\mod{N}$ validates all the pp-sentences valid in $\mod{M}$. From Proposition \ref{lem:atax}, there is a structure $\mod{L}$ elementarily equivalent to $\mod{N}$ and a homomorphism $g: \mod{M} \rightarrow \mod{L}$. Since $T$ is closed under homomorphisms, and $\mod{M}$ is a model of $T$, so is $\mod{L}$. Moreover, since $\mod{N}$ is elementarily equivalent to $\mod{L}$, then $\mod{N}$ is also a model of $T$, proving that $ppT \models_{\mathbf{A}} T$. 

In the case that $\overline{ppT}= \emptyset$, given that $T$ is consistent, for any model of $T$ it holds that $\mod{N}$ validates all the pp-sentences valid in this model (i.e., all possible pp-sentences of $P$, in fact). Then we proceed as above. \end{proof}

Observe that, on the one hand, it is clear that a theory is axiomatized by $pp$-sentences if and only if it is axiomatized by $\wedge$-positive sentences and if and only if it is axiomatized by $\&$-positive sentences. Moreover, since all these formulas are also preserved under weak $\alg{A}$-directed products, we get the following corollary.
\begin{corollary}
	Let $T$ be a consistent theory over $\mathcal{P}$. Then the following are equivalent:
	\begin{itemize}
		\item $T$ is closed under weak direct products and homomorphisms,
		\item $T$ is closed under direct products and homomorphisms,
		\item $T$ is axiomatized by $\wedge$-positive formulas, 
		\item $T$ is axiomatized by $\&$-positive formulas, 
		\item $T$ is axiomatized by pp-formulas.
	\end{itemize}
\end{corollary}

%


Using $\mathbf{A}$-compactness for consequence we can obtain the following corollary of Theorem \ref{t:axiom}. We introduce here the notation of two sentences $\phi$ and $\alpha$ being \emph{1-equivalent} if and only if $Mod _\alg{A}$($\phi$)$=Mod _\alg{A}$($\alpha$). 

\begin{corollary} Let $\mathcal{P}$ be a predicate language and $\phi$ be a consistent sentence. Then, $\phi$ is 1-equivalent to a pp-sentence if and only if $\phi$ is closed under homomorphisms and weak direct products.
\end{corollary}

\begin{proof} One direction is clear. Assume that $\phi$ is a consistent sentence such that Mod($\phi$) are closed under homomorphisms and weak direct products. By Theorem \ref{t:axiom}, the set of consequences of $\phi$ is axiomatized by a set of pp-sentences $\Sigma$. By Theorem \ref{com}, since $\Sigma \models_{\mathbf{A}}\phi$, there is a finite subset $\Sigma_0\subseteq\Sigma$ such that $\Sigma_0 \models_{\mathbf{A}}\phi$. Since pp-sentences are closed under $\wedge$, a pp-sentence $\alpha$ equivalent to the conjunction of all sentences in $\Sigma_0$, is 1-equivalent to $\phi$. \end{proof}

\section{Fuzzy Existential Positive Sets of Axioms}
It is a natural question to ask which kind of sentences axiomatize those classes closed simply under homomorphisms. Using some previous results, we show that these kind of classes are axiomatized by fuzzy existential positive sentences.
\begin{definition}\label{d2:expp} \textbf{\textit{Fuzzy Existential Positive Formula}} Given a predicate language $\mathcal{P}$, and a $\mathcal{P}$-formula $\phi$, it is said that $\phi$ is \emph{existential positive}, if $\phi$ is of the form $(\exists \overline{x})\psi$, where $\psi$ is a quantifier-free formula built from atomic formulas by using only the connectives $\wedge$, $\vee$ and $\&$.
\end{definition} 

Notice that the set of existential formulas is closed under weak disjunction, that is, any finite disjunction of existential formulas in equivalent in MTL$\forall_=$ to
an existential formula (for a reference see p. 281 item (16) of \cite{EsGo01}). Using this fact, the following result of \cite{EsGo01} and the next proposition, it can be proved that every existential positive formula is equivalent to a disjunction of pp-formulas.
\begin{proposition}\label{godo2} \textit{\textbf{Prop. 1.30 and 1.31, \cite{EsGo01}}} Let $\mathcal{P}$ be a predicate language, for every $\mathcal{P}$-formulas $\phi, \psi$ and $\alpha$, the formulas $\phi \& (\psi \vee \alpha) \leftrightarrow ((\phi \& \psi ) \vee (\phi \& \alpha))$, $\phi \wedge (\psi \vee \alpha) \leftrightarrow ((\phi \wedge \psi ) \vee (\phi \& \alpha))$ and $\phi \vee (\psi \wedge \alpha) \leftrightarrow ((\phi \vee \psi ) \wedge (\phi \vee \alpha))$ are \emph{MTL}$\forall_=$ theorems.
\end{proposition}

Now we prove an axiomatization theorem for classes closed under fuzzy homomorphisms.

\begin{theorem}\label{t:axiomex} Let $\mathcal{P}$ be a predicate language and $T$ be a consistent theory. Then, $T$ is closed under homomorphisms if and only if $T$ is axiomatized by a set of existential positive sentences.
\end{theorem}

\begin{proof} $\Leftarrow$ is immediate. For the other direction, consider the set of sentences $epT = \{\varphi: \varphi$ is a disjunction of pp-sentences and $T \models_{\mathbf{A}}\varphi\}$. We aim to prove that $epT \models_{\mathbf{A}} T$. Clearly $epT \neq \emptyset$. Consider a model $\mod{N}$ of $epT$. Notice that, for every pp-sentence $\alpha$, if $||\alpha||_{\mod{N}}<1$, then there is a model $\mod{M}$ of $T$ such that $||\alpha||_{\mod{M}}<1$, otherwise $\alpha \in epT$.

Now expand the language adding a new constant symbol $\overline{a}$, for the coatom of the algebra (the maximum element $a$ such that $a<1$). Let us denote by $\mathcal{P}_a$ this expanded language. Consider the following set of $\mathcal{P}_a$-sentences:
$$\Sigma= \{\alpha \to \overline{a}: \alpha \text{ is a pp-sentence of language } \mathcal{P} \text{ and } ||\alpha||_{\mod{N}}<1\}$$

\noindent We show that $T \cup \Sigma$ has a standard model for the language $\mathcal{P}_a$. By $\alg{A}$-compactness for satisfiability, it is enough to prove that for every finite subset $\Sigma_0 \subseteq \Sigma$, $\Sigma_0=\{\alpha_1 \to \overline{a} \ldots, \alpha_k \to \overline{a}\}$, $T \cup \Sigma_0$ has an standard model for the language $\mathcal{P}_a$.

Assume, searching for a contradiction, that $T \cup \Sigma_0$ has no standard model in the expanded language. That is, in every standard model of $T$ for $\mathcal{P}_a$, for every $1 \leq i \leq k$, the sentence $\alpha_i$ is valid. Observe that every model of $T$ for $\mathcal{P}$ can be expanded to a standard model for the new language $\mathcal{P}_a$ in a natural way (by interpreting the new constant as the corresponding element of the algebra). Consequently, since for every $1 \leq i \leq k$, $\alpha_i$ is a $\mathcal{P}$-sentence, in every model of $T$ for $\mathcal{P}$ the $\alpha_i$'s will be valid. Therefore we will got to a contradiction, we will have that $T \models_{\mathbf{A}}\bigvee_{1 \leq i \leq k} \alpha_i$.

By $\alg{A}$-compactness for satisfiability, there is a standard model of $T \cup \Sigma$, whose reduction $\mod{M}$ to the original language has the property that $\mod{N}$ validates all the pp-sentences valid in $\mod{M}$. From Proposition \ref{lem:atax}, there is a structure $\mod{L}$ elementarily equivalent to $\mod{N}$ and a homomorphism $g: \mod{M} \rightarrow \mod{L}$. Since $T$ is closed under homomorphisms, and $\mod{M}$ is a model of $T$, so is $\mod{L}$. Moreover, since $\mod{N}$ is elementarily equivalent to $\mod{L}$, then $\mod{N}$ is also a model of $T$, proving that $epT \models_{\mathbf{A}} T$. \end{proof}

Following the proof of Theorem \ref{t:axiomex} one can see that it is possible to subtitute in its statement \emph{existential positive sentences} by \emph{existential positive sentences when only the weak conjunction and disjunction occurs}.
\bigskip

Using $\mathbf{A}$-compactness for consequence we can obtain the following corollary of Theorem \ref{t:axiomex}.

\begin{corollary} Let $\mathcal{P}$ be a predicate language and $\phi$ be a consistent sentence. Then, $\phi$ is 1-equivalent to an existential positive sentence if and only if Mod($\phi$) is closed under homomorphisms.
\end{corollary}

\section{Discussion and Future Work}

Can non-classical logic contribute to the analysis of complexity in computer science? We started the paper with the statement of this general question, and in this final section, we would like to comment on how the axiomatization theorem can be regarded as a contribution to provide an answer to it. 

In one of the books of reference in the field \cite{Chang73} model theory is described as algebra+logic. Working in this same framework, and in the line of recent works taking an algebraic approach to valued CSP (see for instance \cite{KolKroRo} and \cite{KroZiv17}), we have presented an algebraic characterization of the preservation of pp-formulas in terms of weak direct products and homomorphisms. Theorem \ref{t:axiom} tells us that there is a good trade-off between algebra and logic in the fuzzy positive-primitive fragment. 

However, the notion of fuzzy homomorphism traditionally used in the fuzzy literature, do not encompass other notions of polymorphism such as weighted or fractional polymorphisms (see for instance \cite{KolKroRo} or \cite{KroZiv17}). It seems natural that each of these notions has been born with a different purpose (preserving tractability, preserving truth degrees, preserving truth only, etc), and further research is needed to study other definitions of homomorphism (for example \cite{HorMoVi17} \cite{De11}, or \cite{DeGaNo16}) and see their pros and contras in different contexts (higher expressivity-non polynomial preservation of complexity, lower expressivity - good computational behaviour, etc). Theorem \ref{t:axiom} also sheds light to the fact that, if we introduce stronger notions of homomorphisms, we will need to redefine pp-formulas, possibly using a language expanded with constant symbols for the elements of the valued structure, in order to maintain the correspondence between algebra and logic. 

The relational structures we have studied are over finite algebras, but we have proven the results both, for finite and for infinite domains, in order to cope with applications on infinite templates. In the classical case, the pp-preservation problem restricted to finite structures was solved by B. Rossman in \cite{Ros08} with some previous results, for instance in \cite{At06} in the context of CSP dualities. It would be interesting to prove the corresponding version in the fuzzy context, especially taking into account the improvements recently introduced in \cite{Ros17} with respect to the bounds on the quantifier-rank of the sentences.

Two other important lines of research are under ongoing work, but their more general range of application requires of a deeper and longer process of study and development. One concerns the behaviour of transformations like homomorphisms and direct products over infinite algebras (where the witnessing condition, as well as compacity, are likely to fail). The notion of homomorphism, as we said before, must be designed depending on its purposes, and understanding the a natural generalization to according which context is already a challenging question. Another ongoing work is that of understanding the behaviour of homomorphisms and products whose result is not a structure of the original algebra, but rather considering the algebras in a variety, and moving between structures over those algebras. This problem partially encompasses the previous one, since we will be bounded, in the majority of cases, to work with infinite algebras (as long as we do some reasoning over an infinite direct product, as for instance we do in the proof of Theorem \ref{lem:atax}).

\section{Acknowledgements}

The research leading to these results has received funding from RecerCaixa 2018 project AppPhil. This project has also received funding from the European Union's Horizon 2020 research and innovation program under the Marie Sklodowska-Curie grant agreement No 689176 (SYSMICS project), and by the projects RASO TIN2015-71799-C2-1-P, CIMBVAL TIN2017-89758-R, and the grant 2017SGR-172 from the Generalitat de Catalunya. 

This work was supported by the grant no. CZ.02.2.69/0.0/0.0/17\_050/0008361 of the Operational programme Research, Development, Education of the Ministry of Education, Youth and Sport of the Czech Republic, co-financed by the European Union, and by the grant GA17-04630S of the Czech Science Foundation.


\begin{thebibliography}{5}


\bibitem{At06}
Atserias, A., Dawar A., Kolaitis, Ph. G.: On preservation under homomorphisms and
unions of conjunctive queries. J. ACM, 53(2):208--237 (2006)

\bibitem{BellZa70}
Bellman, R.E., Zadeh, L.A.:
Decision-making in a fuzzy environment.
Manag. Science, 17:141--164 (1970)


\bibitem{Bul17}
Bulatov, A.:
A Dichotomy Theorem for Nonuniform CSPs. In proceedings of
IEEE 58th Annual Symposium on Foundations of Computer Science (FOCS), 319--330 (2017)


\bibitem{BisMonRo97}
Bistarelli, S., Montanari, U., Rossi, F.: 
Semiring-based constraint satisfaction and optimization.
J. ACM, 44(2):201–236 (1997)

\bibitem{Chang73}
Chang, C.C., Keisler, H.J.:
\newblock Model Theory.
\newblock Elsevier Science Publishers (1973)

\bibitem{CiHa10}
Cintula, P., H\'ajek, P.:
\newblock Triangular norm based predicate fuzzy logics.
\newblock Fuzzy Sets and Systems, 161:311--346 (2010)

\bibitem{CiHaNo11}
Cintula, P., H\'ajek, P., Noguera, C. (eds.):
\newblock Handbook of Mathematical Fuzzy Logic, volume 37,
\newblock Studies in Logic, Mathematical Logic and Foundations (2011)

\bibitem{De11}
Dellunde, P.:
Preserving mappings in fuzzy predicate logics. Journal of Logic and Computation, 22(6):1367-1389 (2011).

\bibitem{De14}
Dellunde, P.:
\newblock Applications of ultraproducts: from compactness to fuzzy elementary classes. Logic Journal of the IGPL 22(1):166-180 (2014)

\bibitem{De18}
Dellunde, P.:
Fuzzy Positive Primitive Formulas. MDAI'18: 156-168 (2018)
 
\bibitem{DeGaNo16}
Dellunde, P., Garc\'ia-Cerda\~{n}a, A., Noguera, C.:
L\" owenheim-Skolem theorems for non-classical first-order algebraizable logics.
Log. J. IGPL 24(3):321-345 (2016) 

\bibitem{DuFarPra93}
Dubois, D., Fargier, H., Prade, H.:
The calculus of fuzzy restrictions as a basis for flexible constraint satisfaction. 
In 2nd IEEE Int. Conf. on Fuzzy Systems. IEEE (1993)

 \bibitem{EsGo01}
Esteva, F., Godo, L.: Monoidal t-norm based logic: towards a logic for left-continuous t-norms.
Fuzzy Sets and Systems, 124:271--288 (2001)


\bibitem{Fa74}
Fagin, R.: 
Generalized First-Order Spectra and Polynomial-Time Recognizable Sets. Complexity of Computation, ed. R. Karp, SIAM-AMS Proceedings, 7:27–41 (1974)

R. Fagin. G
\bibitem{FarLang93}
Fargier H., Lang, J.: 
Uncertainty in constraint satisfaction problems: a probabilistic approach.
In Proceedings of the European Conference on Symbolic and Quantitative Approaches to
Reasoning and Uncertainty, 747:97--104 (1993)

\bibitem{Gei68}
Geiger, D.:
Closed systems of functions and predicates. Pacific Journal of Mathematics, 27(1) (1968)


\bibitem{Ha07}
H{\'a}jek, P.:On witnessed models in fuzzy logic.
Mathematical Logic Quarterly, 53(1), pages 66--77 (2007)

\bibitem{Ha07a}
H{\'a}jek, P.:On witnessed models in fuzzy logic II.
Mathematical Logic Quarterly, 53(6), pages 610--615 (2007)


\bibitem{h}
Hodges, W.: Model Theory. Cambridge (1993)

\bibitem{KoVar08}
Kolaitis, P. G., Vardi, M.Y.: 
A Logical Approach to Constraint Satisfaction.
In Complexity of Constraints - An Overview of Current Research Themes,
LNCS 5250, pages 125--155 (2008)

\bibitem{KolKroRo}
Kolmogorov, V., Krokhin, A., Rolinek, M.:
The Complexity of General-Valued CSPs. In FOCS, pages 1246--1258 (2015)

\bibitem{KroZiv17}
Krokhin, A.A., Zivny, S.:
The Complexity of Valued CSPs. The Constraint Satisfaction Problem 2017: 233--266.


\bibitem{MeRoSchiex06}
Meseguer, P., Rossi, F., Schiex, T.: 
Soft constraints. In Handbook of constraint programming,
chapter 9, pages 281--328 (2006).

\bibitem{MouPra98}
Moura J., Prade, H.:
Logical Analysis of Fuzzy Constraint Satisfaction Problems. 
In 7nd IEEE Int. Conf. on Fuzzy Systems. IEEE (1993)

\bibitem{PiRoVen17}
Pini, M.S., Rossi, F., Venable, K.R.:
Compact Preference Representation via Fuzzy Constraints in Stable Matching Problems. ADT 2017: 333--338.


\bibitem{Prest88}
Prest, M.: Model Theory of Modules. Cambridge (1988)


\bibitem{RoBreWalsh11}
Rossi, F., Brent, K., Walsh, T.: A Short Introduction to Preferences,
Synthesis Lectures on Artificial Intelligence and Machine Learning,
Morgan and Claypool Pub. (2011)


\bibitem{HorMoVi17}
Horc{\'{\i}}k, R., Moraschini, T., Vidal, A.: 
An Algebraic Approach to Valued Constraint Satisfaction.
In 26th {EACSL} Annual Conference on Computer Science Logic, 42:1--42:20 (2017)

\bibitem{Ros08}
Rossman, B.: Homomorphism preservation theorems. J. ACM 55(3): 15:1-15:53 (2008)

\bibitem{Ros17}
Rossman, B.:
An Improved Homomorphism Preservation Theorem From Lower Bounds in Circuit Complexity. ITCS 2017, 27:1-27.

\bibitem{Rutt94}
Ruttkay, Z.: Fuzzy constraint satisfaction. In 3rd IEEE Int. Conf. on Fuzzy Systems. IEEE (1994)

\bibitem{SchiexFarVer}
Schiex, T., Fargier, H., Verfaillie, G.: 
Valued constraint satisfaction problems: Hard and easy problems. 
In Proceedings of the Fourteenth International Joint Conference on Artificial Intelligence
(IJCAI 95) pages 631--639 (1995)

\bibitem{To98}
Torra, V.:
On Considering Constraints of Different Importance in Fuzzy Constraint Satisfaction Problems. 
International Journal of Uncertainty, Fuzziness and Knowledge-Based Systems 6(5):489-502 (1998)



\bibitem{Zh17}
Zhuk, D.: The Proof of {CSP} Dichotomy Conjecture. In proceedings of IEEE 58th Annual Symposium on Foundations of Computer Science (FOCS), pp. 331--342 (2017)



\end{thebibliography}
\end{document}